\documentclass{article}
\usepackage{times,amsmath,amssymb,epsf,rotating}

\newtheorem{lemma}{Lemma}
\newtheorem{theorem}{Theorem}

\def\qed{\hspace*{\fill}{\large$\square$}}

\def\nb{^\vee}
\def\arr{^{\makebox[0pt][l]{$\scriptstyle\nearrow$}\swarrow}}
\def\binomial#1#2{{{#1}\choose{#2}}}

\let\phi\varphi
\let\le\leqslant
\let\ge\geqslant

\begin{document}

\title{The asymptotic complexity of partial sorting\\
  \large\em How to learn large posets by pairwise comparisons}

\author{Jobst Heitzig}

\maketitle

\begin{abstract}
The expected number of pairwise comparisons needed to learn a partial order on
$n$ elements is shown to be at least $n^2/4-o(n^2)$, and an algorithm is given
that needs only $n^2/4+o(n^2)$ comparisons on average. In addition, the optimal
strategy for learning a poset with four elements is presented.

~

\noindent
{\bf Key words:} algorithm, average-case, oracle, poset, sorting.
\end{abstract}

\section{Introduction}

Traditional sorting is about learning a linear order. Its complexity is
often measured by the number of pairwise comparisons a sorting algorithm
needs on average, which is known to be $\Theta(n\log n)$. 
It is a straightforward generalization to ask
for algorithms which learn a {\em partial} order by pairwise comparisons, a task
that could be termed {\em partial sorting}. 
Let us designate the set of all strict partial orders on $n=\{0,1,\dots,n-1\}$ by
$\mathcal P(n)$. This set has $2^{n^2/4+o(n^2)}$ many
elements (cf.~\cite{KR}), 
and each pairwise comparison of elements of $n$ has at most three
possible results. A trivial lower bound for the expected number of comparisons
needed to learn some $P\in\mathcal P$ is therefore $\log_3|\mathcal
O(n)|=\frac{n^2}{4\log_2 3}+o(n^2)$, since in a rooted tree with $\ell$ leaves
in which each node has at most $r$ children, the average leaf-root-distance is
at least $\log_r\ell$.

In this paper, a lower bound of $\frac{n^2}{4}-o(n^2)$ is proved, which is
larger than the above by a factor of $\log_2 3\approx 1.58$. In other words,
any learning algorithm for large posets must expect to compare at least about
half of all pairs. Moreover, it will be shown that there are indeed algorithms 
whose expected running time is just $\frac{n^2}{4}+o(n^2)$. Both results use
the fact that for (very) large $n$, almost all posets have a specific
three-leveled shape.  

To underline the asymptotic nature of the results presented below,
Figure~\ref{fig:strat} shows as a contrast the optimal poset learning strategy
for $n=4$ which has been determined by a recursive computer search. Each
node is a possible state (up to (dual) isomorphisms), and the node's diagram
shows all relations (like in a Hasse diagram) and all incomparabilities
(represented by dotted lines) known in that state. The edges show which states can arise from
which others, where loops indicate dualization. Those states in which there is
only one possible type of comparison are framed with thinner lines, so the
other nodes already determine the actual strategy. Its average running
time is $5.461$ comparisons, compared to $6$ pairs and a trivial lower bound of
$\log_3|\mathcal P(4)|=\log_3 219\approx 4.905$ comparisons. The optimal strategy
for $n=5$ takes $8.744$ comparisons on average, while $\binomial 5 2=10$ and
$\log_3|\mathcal P(5)|=\log_3 4231\approx 7.601$.

\begin{sidewaysfigure}[e]\caption{Optimal learning strategy for 4-element posets}\label{fig:strat}
\begin{center}
\epsfbox{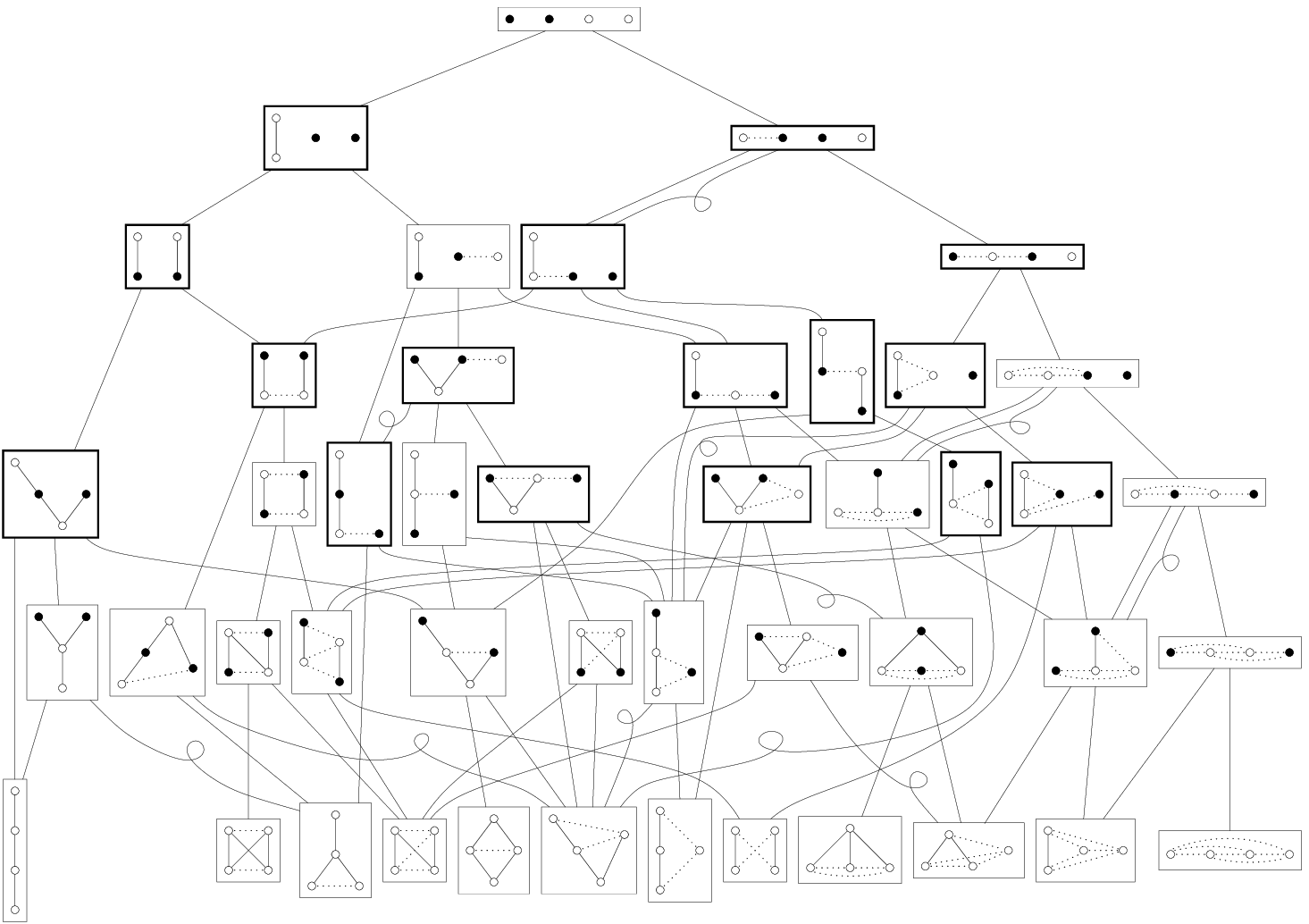}
\end{center}
\end{sidewaysfigure}

\section{The lower bound}

Given $P\in\mathcal P(n)$ and $a,b\in n$, the {\em pairwise comparison} $\{a,b\}$
determines $P|_{\{a,b\}}$, that is, provides the information whether
$a\,P\,b$, $b\,P\,a$, or neither.
Let us define the {\em covering} and {\em anti-covering} relations of
$P$ by
\[ P\nb := P\setminus P^2 \quad\mbox{and}\quad 
   P\arr := \{(x,y)\in n\times n: x\neq y,~Px\subseteq Py,\mbox{~and~} yP\subseteq
   xP\}\setminus P, \]
where $Py=\{x: x\,P\,y\}$ and $xP=\{y: x\,P\,y\}$. 
We consider algorithms which learn a partial order $P\in\mathcal P(n)$ given a 
number $n\ge 1$ and an {\em oracle} for $P$, which is just a subroutine that
performs a pairwise comparison in $P$. The algorithms can learn $P$ only
through oracle calls, each of which is assumed to take constant time.
For any such algorithm $\phi$, let $c_\phi(P)$ be the number of
pairwise comparisons the algorithm needs until it knows $P$. Then
$e_\phi(n):=\sum_{P\in\mathcal P(n)}c_\phi(P)/|\mathcal P(n)|$ is the expected
number of pairwise comparisons for that algorithm.
Finally, let $Q(P):=\big\{\{a,b\}:a\,P\nb\,b$ or $a\,P\arr\,b\big\}$. 
\begin{lemma}
  $c_\phi(P)\ge|Q(P)|$ for all $\phi$ and $P$.
\end{lemma}
{\em Proof.}
Assume that $\phi$ claims to know $P$ but has not compared the pair $\{a,b\}\in
Q(P)$. If $a\,P\nb\,b$, put $P':=P\setminus\{(a,b)\}$, while if $a\,P\arr\,b$,
put $P':=P\cup\{(a,b)\}$. Then $P'$ is a partial order that would erroneously
be recognized as $P$ by $\phi$.\qed 

~

For $R\in\mathcal P(4)$, for example, the average cardinality of $Q(R)$ is
about $4.849$ which is smaller than the trivial lower bound of $4.905$. But
for $R\in\mathcal P(5)$ it is about $7.958$ which improves the trivial
lower bound of $7.601$. 

For the rest of this section, assume that $n$ is a multiple of $4$. 
Let $L(n)$ be the set of all ordered partitions $(A,B,C)$ of $n$ with
$|A|=|C|=n/4$ and $|B|=n/2$. 
Put $\mathcal T(n):=\bigcup_{(A,B,C)\in L(n)}\mathcal T_{ABC}(n)$, where 
$\mathcal T_{ABC}(n)$ is the set of all $P\in\mathcal P(n)$ which fulfil
(i) $x=y$ or $Px\not\subseteq Py$ or $yP\not\subseteq xP$
for all $(x,y)\in A^2\cup B^2\cup C^2$, and
(ii) $aP\cap B\cap Pc\neq\emptyset$ and $A\cap Pb\neq\emptyset\neq bP\cap C$
for all $(a,b,c)\in A\times B\times C$. In particular, these posets consist of
a lower level $A$ of $n/4$ minimal elements, an antichain $B$ of size $n/2$
building the middle level, and an upper level $C$ of $n/4$ maximal elements,
and no $C$-element covers an $A$-element. Moreover, (i) and (ii) imply that
$Q(P)=Q_{ABC}:=(A\times B)\cup(B\times C)$. 

\begin{lemma}\label{lem:3}
  $\displaystyle\frac{|\mathcal T(4m)|}{|\mathcal P(4m)|}=1-o(1/m)$.
\end{lemma}
{\em Proof.} Let $n=4m$. Improving upon the original asymptotics of Kleitman and
Rothschild \cite{KR}, Brightwell, Pr\"omel, and Steger \cite{BPS}
showed that for some $K>1$, $|\mathcal P(n)|=|\mathcal S(n)|\big(1+O(K^{-n})\big)$, where
$\mathcal S(n):=\bigcup_{(A,B,C)\in L(n)}\mathcal S_{ABC}(n)$ and
$\mathcal S_{ABC}(n)$ is the set of all $P\in\mathcal P(n)$ with
$P\nb\subseteq Q_{ABC}$ and $|Px|>1$ for all $x\in B\cup C$. 
On the other hand, it is easy to see that $\mathcal T(n)\subseteq\mathcal S(n)$ and  
$|\mathcal T(n)|/|\mathcal S(n)|=1-o(1/n)$. Hence 
$1-|\mathcal T(n)|/|\mathcal P(n)|=1-\frac{1-o(1/n)}{1+O(K^{-n})}=o(1/n)$.\qed

~

Because $P\in\mathcal T(n)$ implies $|Q(P)|=n^2/4$, it follows 
that $e_\phi(n)$ has a lower bound of $n^2/4-o(n^2)$.
Table~\ref{tab:1} compares $n^2/4$ with $\log_3|\mathcal P(n)|$ for some
small values of $n$ (based on numbers from \cite{HR}).

\begin{table}[e]\caption{Comparison of lower bounds for $e_\phi(n)$ for small $n$}\label{tab:1}
\begin{center}
\begin{tabular}{l|l|l|l|l|l|l|l}
$n$ & 4 & 8 & 12 & 13 & 14 & 15 & 16 \\
$n^2/4$ & 4 & 16 & 36 & 42.25 & 49 & 56.25 & 64 \\
$\log_3|\mathcal P(n)|$ & 4.91 & 18.10 & 36.93 & 42.41 & 48.19 & 54.26 & 58.52
\end{tabular}
\end{center}
\end{table}

\section{A simple algorithm}

Consider the algorithm $\phi_3$ listed in Figure~\ref{fig:alg} which learns a
partial order on $N$. If $N$ is a multiple of
$4$, the strategy of $\phi_3$ first assumes that $P$ is a member of $\mathcal
T(N)$. If the assumption is true, $\phi_3$ will determine the corresponding
level partition $(A,B,C)$ in $o(n^2)$ expected time so that it can afterwards
compare exactly the $n^2/4$ pairs in $Q(P)=Q_{ABC}$. 
In the asymptotically unlikely case that $P\notin\mathcal
T(N)$ it will detect that fact and perform a comparison of all pairs.

Although this is obviously not the best possible strategy, the amount of time
$\phi_3$ ``wastes'' becomes negligible for $N\to\infty$.

\begin{figure}[e]\caption{The asymptotically optimal algorithm $\phi_3$}\label{fig:alg}

\begin{tabbing}
~~~~~~~~~\=input:~~~~~\=oracle $\mathcal C$ for pairwise comparisons in a partial order $P$ on $N$\\
\>output:\>$P$\\[2mm]
\>~~1~~~~~\=put $A=B=C=\emptyset$\\
\>~~2\>find largest $n\le N$ with $4|n$\\[1mm]
\>\>   {\em main loop:}\\
\>~~3\>{\bf for} $k$ from $0$ to $n-1$ {\bf do}\\
\>~~4\>~~~~~\=put $r=|A\cup C|$ and $m=|B|$\\
\>~~5\>\>     assume \=$A\cup C=\{x_1,\cdots,x_r\}$, $B=\{y_1,\cdots,y_m\}$,\\
\>\>\>\>              and $n\setminus\{k\}=\{x_1,\dots,x_{n-1}\}=\{y_1,\dots,y_{n-1}\}$\\[1mm]
\>\>\>        {\em inner loop:}\\
\>~~6\>\>     {\bf for} $i$ from $1$ to $n-1$ {\bf do}\\
\>~~7\>\>~~~~~\={\bf call} $\mathcal C(k,x_i)$\\  
\>~~8\>\>\>     {\bf if}, for some $j<i$, either\\
\>~~9\>\>\>~~~~~\=$x_i\,P\,k\,P\,x_j$, or\\
\>10\>\>\>\>     $x_j\,P\,k\,P\,x_i$, or\\
\>11\>\>\>\>     ($x_i\,P\,k$, $x_j\in A$, but not $x_j\,P\,k$), or\\
\>12\>\>\>\>     ($k\,P\,x_i$, $x_j\in C$, but not $k\,P\,x_j$)\\
\>13\>\>\>\>~~~~~{\bf then} add $k$ to $B$ and continue in main loop\\
\>14\>\>\>     {\bf call} $\mathcal C(k,y_i)$\\  
\>15\>\>\>     {\bf if} $i\le m$ and $k\,P\,y_i$~~~\={\bf then} add $k$ to $A$ and continue in main loop\\
\>16\>\>\>     {\bf if} $i\le m$ and $y_i\,P\,k$\>   {\bf then} add $k$ to $C$ and continue in main loop\\
\>\>\>       {\bf end} {\em (of inner loop)}\\[1mm]
\>17\>\>     {\bf if} $k$ is maximal ($\Leftrightarrow$~$k\,P\,y_i$ for no $i$)~~~\={\bf then} add $k$ to $C$\\
\>18\>\>     {\bf if} $k$ is minimal ($\Leftrightarrow$~$y_i\,P\,k$ for no $i$)\>   {\bf then} add $k$ to $A$\\
\>\>       {\bf end} {\em (of main loop)}\\[1mm]
\>19\>     {\bf for} all $(x,y)\in Q_{ABC}\cup (n\times (N\setminus n))$~~~{\bf call} $\mathcal C(x,y)$\\
\>20\>     {\bf if} the calls so far did not determine $P$ uniquely\\
\>21\>~~~~~{\bf then~~~for} all remaining pairs $(x,y)$~~~{\bf call} $\mathcal C(x,y)$\\
\>22\>     compute and print $P$.
\end{tabbing}
\end{figure}

\begin{theorem}
  $\phi_3$ is an asymptotically optimal poset learning algorithm in the
  sense that $e_{\phi_3}(N)=N^2/4+o(N^2)$.
\end{theorem}
{\em Proof.}
Let $P\in\mathcal P(N)$. 
Because of lines 20--21, $\phi_3$ learns $P$ completely. 

Let $\mathcal U(n):=\bigcup_{(A,B,C)\in L(n)}\mathcal
U_{ABC}(n)\supseteq\mathcal T(n)$, where $\mathcal U_{ABC}(n)\supseteq\mathcal
T_{ABC}(n)$ is the set of all $P\in\mathcal P(n)$ with $P\nb\subseteq Q_{ABC}$. 
We may assume that $P|_n\in\mathcal U_{A_0B_0C_0}(n)$ for some
$(A_0,B_0,C_0)\in L(n)$, since by Lemma~\ref{lem:3}, $P|_n\in\mathcal U(n)$ is true with a
probability converging to 1 as $N\to\infty$. 
Note that $\alpha_n:=1-|\mathcal T(n)|/|\mathcal U(n)|=o(1/n)$ is an upper bound
for the probability that at some point in $\phi_3$, either $A\not\subseteq
A_0$, $B\not\subseteq B_0$, or $C\not\subseteq C_0$.

Conditional to $P|_n\in\mathcal U_{A_0B_0C_0}(n)$, the event $x\,P\,y$ has
probability $\frac 1 2$ independently for all $(x,y)\in Q_{A_0B_0C_0}$. 
Hence one can estimate the expected number of
pairwise comparisons in iteration $k$ of the main loop as follows.

(i) Assume that $k\in B_0$. 
For $j:=1$ and $2\le i\le r$, the disjunction in lines 9--12 is violated
with probability at most $\frac 1 2+\alpha_n$. 
Hence iteration $k$ takes an expected number of at most 
\[ 2\left[\sum_{i=2}^r i(\textstyle\frac 1 2+\alpha_n)^{i-2}+(n-1)(\textstyle\frac 1 2+\alpha_n)^{r-1}\right]
   <2\frac{(\textstyle\frac 1 2-\alpha_n)^{-2}+n(\textstyle\frac 1 2+\alpha_n)^r}{\textstyle\frac 1 2+\alpha_n} \]
pairwise comparisons in this case.

(ii) Assume that, on the other hand, $k\in A_0\cup C_0$.
For $1\le i\le m$, the probability that both the conditions of lines 15--16 are
violated is at most $\frac 1 2+\alpha_n$ so that in this case iteration $k$ takes an
expected number of at most  
\[ 2\left[\sum_{i=1}^m i(\textstyle\frac 1 2+\alpha_n)^{i-1}+(n-1)(\textstyle\frac 1 2+\alpha_n)^m\right]
   <2\frac{(\textstyle\frac 1 2-\alpha_n)^{-2}+n(\textstyle\frac 1 2+\alpha_n)^m}{\textstyle\frac 1 2+\alpha_n} \]

We have seen that, given $r$ and $m$, iteration $k$ takes in both cases an expected number of
at most $(16+4n(\frac 1 2+\alpha_n)^a)\beta_n$ pairwise comparisons,
where $a:=\min\{r,m\}$ and $\beta_n\to 1$.
At the beginning of iteration $k$, for $0\le a\le\frac k 2$, the probability
that $r=a$ and $m=k-a$ is at most 
\[ \frac{\binomial{n/2}{a}\binomial{n/2}{k-a}}{\binomial{n}{k}}+\alpha_n
   =\frac{\binomial{k}{a}\binomial{n-k}{n/2-a}}{\binomial{n}{n/2}}+\alpha_n
   <\frac{2^{n-k}}{\binomial{n}{n/2}}\binomial{k}{a}+\alpha_n, \]
and so is the probability that $m=a$ and $r=k-a$. In contrast, the probability
that $r+m\neq k$ is at most $\alpha_n$. In all, iteration $k$ takes an expected
number of at most
\begin{eqnarray*}
  \lefteqn{2\sum_{a=0}^{\lfloor k/2\rfloor}
     \left(\frac{2^{n-k}}{\binomial{n}{n/2}}\binomial{k}{a}+\alpha_n\right)
     (16+4n(\textstyle\frac 1 2+\alpha)^a)\beta(n)+\alpha_n\binomial n 2}\\
   &\le& O(n)\frac{2^{n-k}}{\binomial{n}{n/2}}\sum_{a=0}^k\binomial k a
     (\textstyle\frac 1 2+\alpha_n)^a 1^{k-a}+o(n)\\
   &=& O(n)\frac{2^n}{\binomial{n}{n/2}}
     (\textstyle\frac 3 4+\frac{\alpha_n}{2})^k+o(n)
   = O(n^{3/2})(\textstyle\frac 3 4+\frac{\alpha_n}{2})^k
\end{eqnarray*}
pairwise comparisons, so that the total expected number of comparisons in
lines 1--18 is $O(n^{3/2})$. If $P\in\mathcal T(N)$ then $P$ is uniquely
determined in line 20, hence the expected number of comparisons in
lines 19--21 is $N^2/4+o(N^2)$, proving the theorem.\qed

\end{document}